\documentclass[11pt, a4paper]{amsart}

\usepackage{amsfonts,amsmath,amssymb,amscd}

\usepackage[all]{xy}

\newcommand{\lbtimess}{\mathbin{\raisebox{0.2ex}{\makebox[1em][l]{${\scriptstyle\blacktriangleright\mathrel{\mkern-4mu}<}$}}}}

\newcommand{\ldtimess}{\mathbin{\raisebox{0.2ex}{\makebox[1em][l]{${\scriptstyle\rlap{$\scriptstyle{\gtrdot}$}{\vartriangleright}\mathrel{\mkern-4mu}<}$}}}}

\newtheorem{theorem}{Theorem}[section]
\newtheorem{lemma}[theorem]{Lemma}
\newtheorem{prop}[theorem]{Proposition}

\theoremstyle{definition}
\newtheorem{definition}[theorem]{Definition}
\newtheorem{rem}[theorem]{Remark}

\newtheorem{example}[theorem]{Example}

\usepackage{color}

\numberwithin{equation}{section}

\title[Super algebraic groups]
{Hopf algebraic techniques applied to\\ super algebraic groups}

\author{Akira Masuoka}
\address{Akira Masuoka: 
Institute of Mathematics, 
University of Tsukuba, 
Ibaraki 305-8571, Japan}
\email{akira@math.tsukuba.ac.jp}

\begin{document}

\begin{abstract}
Reproducing my talk at Algebra Symposium held at Hiroshima University,
August 26--29, 2013, I review 
recent results on super algebraic groups, emphasizing results 
obtained by myself and my coauthors   
using Hopf algebraic techniques. The results are all basic, and I intend to
make this report into a somewhat informal introduction to the subject. 
\end{abstract}

\maketitle



\section{Motivation}\label{sec:motivation}

We work over $\Bbbk$ which we suppose to be a field (with a very few exceptions).
Thus, vector spaces, tensor products $\otimes$, 
any kind of algebras and so on are supposed to be over the field $\Bbbk$, unless otherwise stated. 

Let us start with the following famous theorem due to Deligne, which is here formulated rather informally.

\begin{theorem}[Deligne \cite{D}]\label{thm:Deligine}
Suppose that $\Bbbk$ is an algebraically closed field of characteristic zero. Then any
rigid symmetric $\Bbbk$-linear abelian tensor category that satisfies some mild, algebraic assumption
is realized as the category of finite-dimensional super modules over some super algebraic group. 
\end{theorem}

The definition of super algebraic groups is quite simple, as will be seen below. But, when I encountered
this theorem around 2003, little seemed to be known about them, compared with super Lie groups which
have a longer history of study founded by Kostant \cite{Kostant}, Koszul \cite{Koszul} and others in the 1970's.  
Indeed, there was not yet proven even the one-to-one correspondence between 
the closed normal super subgroups of a given super algebraic group 
and its quotient super groups; see Section \ref{subsec:super_scheme} below. 
That is why I became, though slowly, to study the subject. 
I use Hopf algebraic techniques, following Hochschild and Takeuchi
who studied algebraic groups by using those techniques in the 1970's; see \cite{H, T1, T2, T3, T4}.  

I am going to review 
recent results on super algebraic groups, emphasizing results 
obtained by myself and my coauthors. 
For my knowledge of the subject I owe very much Alexandr N.~Zubkov, one of my 
coauthors. See \cite{Carmeli-Caston-Fioresi, DM, V}
for modern treatment of wider topics on super geometry. 

\section{What is an affine/algebraic group?}\label{sec:affine_group}

Before going into the super world let us recall what happened in the classical non-super situation.

\subsection{Geometric vs. functorial viewpoints}\label{subsec:G_vs_F}
What is a \emph{space} over $\Bbbk$? It is defined from \emph{geometric viewpoint} to be a locally 
ringed space $(X, \mathcal{O}_X)$, 
where $\mathcal{O}_X$ is a sheaf of commutative algebras over $\Bbbk$. 
From \emph{functorial viewpoint} it is defined to be
a set-valued functor defined on the category $(\mathsf{Comm~Algebras})$ of commutative algebras over $\Bbbk$.
The notion of \emph{schemes} is defined separately from each viewpoint. 
The \emph{Comparison Theorem}\footnote{The theorem is indeed 
proved when $\Bbbk$ is an arbitrary commutative ring.} \cite[I, \S 1, 4.4]{DG}
states that
there exists a category equivalence between the schemes defined from geometric viewpoint and
the schemes defined from functorial viewpoint.

Here we recall how \emph{schemes} or \emph{sheaves} are defined from functorial viewpoint. 
Let $X$ be a set-valued functor $X$ defined on $(\mathsf{Comm~Algebras})$.
We say that $X$ is \emph{affine} if it is representable. A \emph{scheme} (over $\Bbbk$) is a 
local functor defined on $(\mathsf{Comm~Algebras})$ which is covered by open affine subfunctors.  
A map $R \to T$ of commutative algebras is called an \emph{fpqc}
 (resp., \emph{fppf}) \emph{covering},\
if (i)~$T$ is faithfully flat over $R$ (resp., if (i) and (ii)~$T$ is finitely presented as an $R$-algebra). 
Such a map gives a natural exact diagram
\begin{equation}\label{eq:fpqf_covering}
R \to T \rightrightarrows T\otimes_R T. 
\end{equation}
We say that $X$ is a \emph{dur sheaf} (resp., \emph{sheaf}) (over $\Bbbk$) if it preserves 
finite direct products and any exact diagram as given above. Obviously, a dur sheaf 
is a sheaf. 
It is known that a scheme is a dur sheaf; see \cite[III, \S 1, 3.3]{DG}. 

\subsection{Affine groups}\label{subsec:affine_group}
It is convenient and even natural to define affine/alge-\\ braic groups 
from functorial viewpoint. An \emph{affine group} is
a representable functor $G : (\mathsf{Comm~Algebras})\to (\mathsf{Groups})$ with values in the category
of groups. The representing algebra $\mathcal{O}(G)$ has uniquely structure maps of a Hopf algebra 
\[ \Delta : \mathcal{O}(G)\to \mathcal{O}(G)\otimes \mathcal{O}(G),\quad \varepsilon : \mathcal{O}(G)\to \Bbbk,
\quad S : \mathcal{O}(G)\to \mathcal{O}(G), \]
which give the product, the unit and the inverse of $G$, respectively. Thus we have a category anti-isomorphism
\begin{equation}\label{eq:anti-iso}
 (\mathsf{Affine~Groups}) \simeq (\mathsf{Comm~Hopf~Algebras})  
\end{equation} 
between the affine groups and the commutative Hopf algebras. The affine group $G$ which corresponds to
a commutative Hopf algebra is denoted by $\mathrm{Sp}\, A$; this should be once distinguished from the
prime spectrum $\mathrm{Spec}\, A$, though the two notions are eventually equivalent. 

By definition a \emph{closed subgroup} of an affine group $G$ is an affine group $H$
which is represented by a quotient Hopf
algebra of $\mathcal{O}(G)$; it is said to be \emph{normal} if each subgroup $H(R)$ is normal in $G(R)$,
where $R \in (\mathsf{Comm~Algebras})$.  
A \emph{quotient group} of $G$ is an affine group which is represented by
a Hopf subalgebra of $\mathcal{O}(G)$.
This last definition is justified by the 
fact\footnote{We have a very simple, purely Hopf-algebraic proof of this
fact (see \cite{MW}), from which 
the readers may hopefully see that Hopf algebraic techniques are effective.} 
that \emph{a commutative Hopf algebra $A$
is faithfully flat over every Hopf subalgebra $B$}, since it  
implies that $\mathrm{Sp}\, A \to \mathrm{Sp}\, B$
is an epimorphism of dur sheaves. One can prove a natural one-to-one
correspondence between the closed normal subgroups of a given affine group $G$ and the 
quotient groups of $G$. 

A \emph{representation} of an affine group $G$ (or a \emph{left} $G$-\emph{module}) is a morphism of group-valued
functors $\phi : G \to \mathrm{GL}_V$ to the general linear group $\mathrm{GL}_V$ on some vector space $V$;
$\mathrm{GL}_V$ associates to $R \in (\mathsf{Comm~Algebras})$ the group $\mathrm{Aut}_R(V\otimes R)$ of
all $R$-linear automorphism on $V \otimes R$. Since $G$ is representable, such a $\phi$ is uniquely determined
by the image $\phi(\mathrm{id}) : V \otimes \mathcal{O}(G)\overset{\simeq}{\longrightarrow}  
V \otimes \mathcal{O}(G)$ of the identity map on $\mathcal{O}(G)$, or 
by its restriction $\rho := \phi(\mathrm{id})|_V: V
\to V \otimes \mathcal{O}(G)$ to $V = V \otimes \Bbbk$. The requirement that $\phi$ should preserve 
the group structure
is equivalent to that $\rho$ is coassociative and counital, or $(V, \rho)$ is a \emph{right} $\mathcal{O}(G)$-\emph{comodule}. 
In summary, a left $G$-module is the same as a right $\mathcal{O}(G)$-comodule. 

\subsection{Algebraic groups}\label{subsec:algebraic_group}
An affine group $G$ is called an \emph{algebraic group} if $\mathcal{O}(G)$ is finitely generated. A finitely
generated commutative Hopf algebra is called an \emph{affine Hopf algebra}.  
Therefore, the
category anti-isomorphism \eqref{eq:anti-iso} restricts to 
\begin{equation*}
 (\mathsf{Algebraic~Groups}) \simeq (\mathsf{Affine~Hopf~Algebras}).  
\end{equation*} 
Since every 
commutative Hopf algebra is a directed union of finitely generated Hopf subalgebras, every affine
group is a projective limit of algebraic groups. 

Given an algebraic group $G$, then functor points $G(\overline{\Bbbk})$ in the algebraic closure 
$\overline{\Bbbk}$ of $\Bbbk$ form a linear algebraic group over $\overline{\Bbbk}$. If $\Bbbk$ is
an algebraically closed field of characteristic zero, $G \mapsto G(\Bbbk)$ gives an equivalence 
from $(\mathsf{Algebraic~Groups})$ to the category of linear algebraic groups. 

Let $G$ be an algebraic group over a field $\Bbbk$, and let $H$ be a closed
subgroup of $G$. Let $G/H$ denote the functor defined on $(\mathsf{Comm~Algebras})$ which
associates to each $R$ the set $G(R)/H(R)$ of left cosets. We have uniquely a \emph{sheafification}
$G\tilde{/} H$ of $G/H$, that is, a sheaf given a 
functor morphism from $G/H$ that has the obvious universality. We have a natural epimorphism 
$G \to G\tilde{/} H$ of sheaves. Here is a well-known theorem;
see \cite[Part~I, 5.6, (8)]{Jantzen}. 

\begin{theorem}\label{thm:G/H}
$G\tilde{/} H$ is a Noetherian scheme such that $G \to G\tilde{/}H$ is affine and faithfully flat. 
\end{theorem}

The result applied to the opposite algebraic groups $G^{op} \supset H^{op}$ shows an analogous result 
for the sheafification $H \tilde{\backslash} G$ of 
the functor $R \mapsto H(R) \backslash G(R)$ giving right cosets.

\subsection{Hyperalgebras}
Takeuchi studied algebraic groups via characteristic-free approach using Hopf
algebras; see \cite{T1, T2, T3, T4} for example. 
Compared with commutative Hopf algebras, cocommutative Hopf algebras are much more
tractable. Suppose that we are given an algebraic group $G$, and let $A = \mathcal{O}(G)$. 
Takeuchi's main idea is to study the associated cocommutative Hopf algebra $\mathrm{hy}(G)$,
which is called the \emph{hyperalgebra} of $G$. By definition this consists of those
elements in the dual vector space $A^*$ of $A$ which annihilate some powers $(A^+)^n$,
$0 < n \in \mathbb{Z}$, of the augmentation ideal $A^+ := \mathrm{Ker}(\varepsilon : A \to \Bbbk)$. 
Thus we have
\begin{equation}\label{eq:hyperalgebra}
\mathrm{hy}(G)= \bigcup_{n>0} (A/(A^+)^n)^*. 
\end{equation}
This is indeed a cocommutative Hopf algebra which is \emph{irreducible} as a coalgebra
in the sense that 
the trivial $\mathrm{hy}(G)$-comodule $\Bbbk$ is the unique simple $\mathrm{hy}(G)$-comodule.
If the characteristic $\mathrm{char}\, \Bbbk$ is zero, then $\mathrm{hy}(G)$ coincides with the universal 
envelope $U(\mathrm{Lie}(G))$ of the Lie algebra $\mathrm{Lie}(G)$ of $G$. 

A \emph{hyperalgebra} is a synonym of an irreducible cocommutative Hopf algebra, and it
may be regarded as a generalized object of Lie algebras\footnote{I hear that it used to be called
a \emph{hyper-Lie algebra}, before Takeuchi removed  $``$Lie" from the name.}.  
Takeuchi proved that $\mathrm{hy}(G)$ reflects various properties of $G$, even better than $\mathrm{Lie}(G)$ does
in some situations. 
Recently, $\mathrm{hy}(G)$ is often denoted alternatively by $\mathrm{Dist}(G)$,
called the \emph{distribution algebra} of $G$. 
I wish to use $\mathrm{hy}(G)$ in honor of Takeuchi's contributions. 

\subsection{The dual coalgebra}\label{subsec:dual_coalgebra}
Let $A$ be an algebra. Given an ideal $I \subset A$ which is \emph{cofinite} in the sense
$\dim A/I < \infty$, the dual space $(A/I)^*$ of $A/I$ is naturally a coalgebra. Therefore, the directed union
\begin{equation}\label{eq:dual_coalgebra}
A^{\circ}:= \bigcup_{I}\, (A/I)^*\, \subset A^*, 
\end{equation}
where $I$ runs over all cofinite ideals of $A$, is a coalgebra, which is called the \emph{dual coalgebra} of $A$. 
This coincides with the \emph{coefficient space} of all finite-dimensional representations, 
$\pi : A \to \mathrm{End}(V)$, of $A$; it is by definition the union $\bigcup_{\pi}\mathrm{Im}(\pi^*)$ in $A^*$. 
If $A$ is a Hopf algebra, then $A^{\circ}$ is a Hopf algebra. 

Let $V$ be a vector space possibly of infinite dimension. Given a right $A^{\circ}$-comodule structure 
$\rho : V \to V \otimes A^{\circ}$ on $V$, we have a locally finite left $A$-module structure given by
\[ a\, v := \sum_i f_i(a)\, v_i, \quad a\in A,\ v \in V,\]
where $\rho(v)= \sum_i v_i\otimes f_i$.  This gives rise to a bijection from the set of all
right $A^{\circ}$-comodule structures on $V$ to the set of all locally finite left $A$-module structures on $V$.

\section{Invitation to the super world}\label{sec:super_world}
In what follows until end of this report we assume $\mathrm{char}\, \Bbbk \ne 2$.

\subsection{Super vector spaces}
Let $\mathbb{Z}_2=\{ 0, 1 \}$ is the finite group of order 2. 
$``$Super" is a synonym of $``\mathbb{Z}_2$-graded". So, a super vector space is a vector space
$V = V_0 \oplus V_1$ decomposed into subspaces $V_i$ indexed by $i = 0, 1$. The component $V_i$
and its elements are said to be \emph{even} or \emph{odd}, respectively, if $i =0$ or if $i=1$. We say that
$V$ is \emph{purely even} or \emph{odd}, respectively, if $V =V_0$ or if $V=V_1$. 
The super vector spaces $V, W, \dots$ form a tensor category $(\mathsf{Super~Vec~Spaces})$;
the tensor product is that of vector spaces $V \otimes W$ which is $\mathbb{Z}_2$-graded by
the total degree, 
\[ (V \otimes W)_k := \bigoplus_{i+j=k}V_i \otimes W_j,\quad k=0,1. \] 
The unit object is $\Bbbk$ which is purely even. This tensor category is symmetric with respect to the
symmetry $c_{V,W} : V \otimes W \overset{\simeq}{\longrightarrow} W \otimes V$
given by 
\begin{equation*}
c_{V,W}(v \otimes w) = (-1)^{|v|\, |w|}w\otimes v 
= 
\begin{cases}
-w\otimes v &\text{if}\ |v| = |w| = 1\\
w \otimes v & \text{otherwise}.
\end{cases} 
\end{equation*}
Here $v$, $w$ are supposed, as our convention,
to be homogeneous elements with degrees $|v|$, $|w|$. The assumption
$\mathrm{char}\, \Bbbk \ne 2$ ensures that this is different from the obvious symmetry $v \otimes w
\mapsto w \otimes v$. 
The symmetry defined above is called the
\emph{super symmetry}, and will be depicted by 
$$
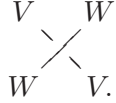
\begin{xy}
(0,0)*+{V}="a";
(10,0)*+{W}="b",
(0,-10)*+{W} ="bb",
(10,-10)*+{V.}="aa",
{"a" \SelectTips{cm}{} \ar @{-}|!{"b";"bb"}\hole "aa"},
{"b" \SelectTips{cm}{} \ar @{-} "bb"},
\end{xy}
$$

\subsection{Super objects}
Algebraic systems, such as algebra or Hopf (or Lie) algebra, defined in the tensor category $(\mathsf{Vec~Spaces})$
of vector spaces 
with the obvious symmetry are generalized by algebraic systems in $(\mathsf{Super~Vec~Spaces})$, which are
called with the prefix $``$super", so as super algebra or super Hopf (or Lie) algebra\footnote{This usage of
terms is not necessarily in fashion, but it seems more acceptable for non-specialists.}. An essential difference
appears when the symmetries are concerned in the argument. For example,
a \emph{super algebra} is just a $\mathbb{Z}_2$-graded algebra, which is not concerned with the super symmetry. 
But the tensor product $A \otimes B$ of super algebras $A$, $B$ involves the super symmetry, 
and its product is given by
\[ (a \otimes b)(c \otimes d) = (-1)^{|b|\, |c|}(ac \otimes bd). \] 
This is depicted by 
$$
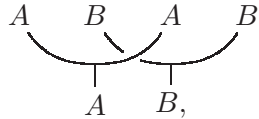
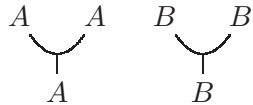
\begin{xy}
(0,0)   *+{A}          ="a",
(20,0)  *+{A}          ="b",
(10,-12)*+{A}          ="c",
(10,0)   *+{B}          ="aa",
(30,0)  *+{B}          ="bb",
(20,-12)*+{B,}          ="cc",
(10,-6.5) *{}="cen",
(20,-6.5) *{}="cent",
(15.5,-5)*{}="hol",
\ar @/_/@{-} "a";"cen"
\ar @/_/@{-} "cen";"b"
\ar @/_/@{-} "aa";"cent"|!{"hol";"b"}\hole
\ar @/_/@{-} "cent";"bb"
\ar @{-} "cen";"c"
\ar @{-} "cent";"cc"
\end{xy}
$$
where
$$\begin{xy}
(0,0)*+{A}="a";
(10,0)*+{A}="b",
**\crv{(5,-10)},
(5,-5)*{} ="cen",
(5,-10)*+{A}="c",
{"cen" \SelectTips{cm}{} \ar @{-} "c"},
\end{xy}
\quad
\begin{xy}
(0,0)*+{B}="a";
(10,0)*+{B}="b",
**\crv{(5,-10)},
(5,-5)*{} ="cen",
(5,-10)*+{B}="c",
{"cen" \SelectTips{cm}{} \ar @{-} "c"},
\end{xy}
$$
represent the products on $A$, $B$, respectively. 
To emphasize this situation, we will write $A\, \underline{\otimes}\, B$ for $A \otimes B$. 

For a super bi- or Hopf algebra $A$,  the coproduct 
$\Delta : A \to  A\, \underline{\otimes}\, A$ is required to be a $\mathbb{Z}_2$-graded algebra map. 

A \emph{super Lie algebra} is a super vector space $\mathfrak{g} = \mathfrak{g}_0 \oplus \mathfrak{g}_1$
given a $\mathbb{Z}_2$-graded linear map 
$[\ , \ ] : \mathfrak{g} \otimes \mathfrak{g} \to \mathfrak{g}$, called \emph{bracket}, which satisfies 
$$[\ , \ ] \circ (\mathrm{id}_{\mathfrak{g} \otimes \mathfrak{g}} + c_{\mathfrak{g},\mathfrak{g}})= 0,\quad
[[\ , \ ],\ ] \circ (\mathrm{id}_{\mathfrak{g}\otimes \mathfrak{g}\otimes \mathfrak{g}} + 
c_{\mathfrak{g}, \mathfrak{g}\otimes \mathfrak{g}} + c_{\mathfrak{g} \otimes \mathfrak{g}, \mathfrak{g}}) = 0.$$
Note that $\mathfrak{g}_0$ is then an ordinary Lie algebra. 

Ordinary objects such as Hopf or Lie algebra are regarded as purely even super objects,
such as purely even super Hopf or Lie algebra. 

A super algebra $A$ is said to be \emph{super-commutative} if the product $A \otimes A \to A$ is invariant,
composed with $c_{A,A}$, or more explicitly, if we have
\[ ab = ba\ \text{if}\ a\ \text{or}\ b\ \text{is even, and}\ ab = - ba\ \text{if}\ a, b\ \text{are both odd}. \]
The last condition is equivalent to that $a^2=0$ if $a$ is odd. Dually, one defines the 
\emph{super-cocommutativity} for super coalgebras. 

A super-commutative super (Hopf) algebra will be called a \emph{super commutative (Hopf) algebra},
regarded as a commutative (Hopf) algebra object in $(\mathsf{Super~Vec~Spaces})$. 

\begin{example}\label{ex:exterior}
Let $V$ be a vector space. The exterior algebra $\wedge (V)= \bigoplus_{n=0}^{\infty}\wedge^n(V)$ is graded by
$\mathbb{N}=\{ 0,1, \dots \}$, and hence is graded by $\mathbb{Z}_2$, so as $(\wedge (V))_i = 
\bigoplus_{m=0}^{\infty}\wedge^{2m+i}(V)$, $i=0,1$. This turns uniquely into a super Hopf algebra in which
each element $v \in V$ is (odd) primitive, i.e. $\Delta(v) = 1\otimes v + v \otimes 1$.
This super Hopf algebra is super-commutative and super-cocommutative. The pairings
$\langle\, ,\, \rangle:\wedge^n(V^*)\times \wedge^n(V)\to\Bbbk$, $n=0,1,\dots$, given by
\begin{equation*}
\langle f_1 \wedge \dots \wedge f_n, \, v_1 \wedge \dots \wedge v_n \rangle = 
\sum_{ \sigma \in \mathfrak{S}_n }(\operatorname{sgn} \sigma) f_1(v_{ \sigma (1) })\,  
\dots f_n(v_{ \sigma (n) }), \quad f_i \in V^*, \ v_i \in V
\end{equation*}
are summarized to $\langle\, ,\, \rangle:\wedge(V^*)\times \wedge(V)\to\Bbbk$, which induces
an isomorphism of super Hopf algebras $\wedge(V^*)\overset{\simeq}{\longrightarrow} \wedge(V)^*$ if
$\dim V< \infty$. 
\end{example}

\subsection{Bosonization}\label{subsec:bosonization}
A super vector space $V$ is identified with a module over
the group algebra $\Bbbk \mathbb{Z}_2$; the generator of $\mathbb{Z}_2$ acts on $V_i$ by scalar
multiplication by $(-1)^i$. A super algebra $A$ is identified with an algebra on which $\mathbb{Z}_2$
acts as algebra automorphisms, so that we have the algebra $\mathbb{Z}_2 \ltimes A$ 
of semi-direct (or smash) product. A \emph{super $A$-module} is an $A$-module object
in $(\mathsf{Super~Vec~Spaces})$; it is identified with an ordinary module over $\mathbb{Z}_2\ltimes A$. 
Given a super right $A$-module $M$ and a super left $A$-module $N$, we have
\[ M \otimes_{\mathbb{Z}_2\ltimes A} N = (M\otimes_A N)_0,\quad M \otimes_{\mathbb{Z}_2\ltimes A} N[1] =
(M\otimes_A N)_1, \]
where $N[1]$ denotes the degree shift of $N$ so that $N[1]_0 =N_1$,\ $N[1]_1 = N_0$. 
This together with the fact that $\mathbb{Z}_2\ltimes A$ is faithfully flat over $A$ proves the following.

\begin{lemma}\label{lem:faithfully_flat}
For $M$ as above the following are equivalent:
\begin{itemize}
\item[(1)] $M$ is (faithfully) flat, regarded as an ordinary right $A$-module;
\item[(2)] $M$ is (faithfully) flat as a right $\mathbb{Z}_2\ltimes A$-module;
\item[(3)] The functor $M \otimes_A$ defined on the category of super left $A$-modules, 
which associates $M \otimes_A N$ to each $N$, is (faithfully) exact.
\end{itemize}
\end{lemma}
An analogous result for super \emph{left} $A$-modules holds true. Note that if $A$ is super-commutative, 
super left $A$-modules and super right $A$-modules are naturally identified. 

We remark that a super (left or right) 
$A$-module is projective in the category of super $A$-modules 
(or equivalently, of $\mathbb{Z}_2\ltimes A$-modules)
if and only if it is projective in the category of ordinary $A$-modules; this holds since the ring extension
$\mathbb{Z}_2\ltimes A\supset A$ is separable. 

Let $C$ be a super coalgebra. A \emph{super $C$-comodule} is a $C$-comodule object in 
$(\mathsf{Super~Vec~Spaces})$; it is identified with an ordinary comodule over the coalgebra
of smash coproduct
\[ \mathbb{Z}_2\lbtimess C. \]
This equals $\Bbbk\mathbb{Z}_2 \otimes C$ as a vector space, and has as its counit the tensor product 
$\varepsilon\otimes \varepsilon$ of the counits. The coproduct 
$ \Delta :  \mathbb{Z}_2\lbtimess C \to (\mathbb{Z}_2\lbtimess C) \otimes (\mathbb{Z}_2\lbtimess C) $
is the left $\Bbbk\mathbb{Z}_2$-module map defined by
\begin{equation}\label{eq:smash_coproduct}
\Delta (1 \otimes c ) = \sum_{(c)} \, (1 \otimes c_{(1)})\otimes (|c_{(1)}| \otimes c_{(2)}),\quad c \in C, 
\end{equation}
where $\Delta(c) = \sum_{(c)} c_{(1)} \otimes c_{(2)}$ denotes the coproduct on $C$. A dual result of
Lemma \ref{lem:faithfully_flat} holds; see \cite[Proposition~2.3]{MZ}. 

Suppose that $A$ is a super Hopf algebra. The smash product and coproduct structures make 
$\Bbbk \mathbb{Z}_2 \otimes A$ into an ordinary Hopf algebra, 
\[ \mathbb{Z}_2 \ldtimess A, \]
which is called the \emph{bosonization} of $A$; this construction was given by Radford
\cite{Radford} in a generalized situation. The construction gives us a useful technique to
derive results on super Hopf algebras from known results on ordinary Hopf algebras;
see \cite[Section 10]{MZ} for example.

\section{Super affine/algebraic groups}

Let us start to discuss these objects of our concern. 
\subsection{Definitions}
Once given the definition of affine/algebraic groups, it is quite easy to define their super analogues. One
has only to replace $(\mathsf{Com~Algebras})$ with the category $(\mathsf{Super~Com~Algebras})$  
of super commutative algebras. Thus, a \emph{super affine group} is a representable functor
$G : (\mathsf{Super~Com~Algebras})\to (\mathsf{Groups})$; it is uniquely represented by
a super commutative Hopf algebra, which we denote by $\mathcal{O}(G)$. Such a $G$ is called a \emph{super
algebraic group} if $\mathcal{O}(G)$ is \emph{affine}, i.e. finitely generated (and super-commutative). 
We have thus a category anti-isomorphism between the super affine groups and the super commutative
Hopf algebras, which restricts to a category anti-isomorphism between the super algebraic groups
and the super affine Hopf algebras. The super affine group which corresponds to a super commutative
Hopf algebra $A$ is denote by $\mathrm{SSp}\, A$; it associates to $R \in (\mathsf{Super~Com~Algebras})$
the group of all super algebra maps $A \to R$. A \emph{closed (normal) super subgroup} of a super affine group $G$
is a super affine group $H$ which is represented by a quotient super Hopf algebra of $\mathcal{O}(G)$ (so that
each $H(R)$ is normal in $G(R)$, where $R \in (\mathsf{Super~Com~Algebras})$).   
Just as in the non-super situation, every super affine
group is a projective limit of super algebraic groups. 
  
By restriction of the domain every functor $G : (\mathsf{Super~Com~Algebras})\to (\mathsf{Groups})$ gives 
rise to a functor $(\mathsf{Com~Algebras})\to (\mathsf{Groups})$, which we denote by $G_{ev}$. 
Suppose that $G= \mathrm{SSp}\, A$ is a super affine/algebraic group. Then $G_{ev}$ is an
affine/algebraic group, being represented by 
\[ \overline{A}:= A/(A_1). \]
This is the (largest) quotient purely even super algebra of $A$ divided by the ideal generated by $A_1$, and
is indeed a quotient super Hopf algebra. 
Therefore, $G_{ev}$ can be identified with the closed super subgroup
of $G$ given by $R \mapsto G(R_0)$. 
We will say that $G_{ev}$ is \emph{associated with} $G$. 
This $G_{ev}$ will be seen to play an important role when we study $G$. 

\subsection{Super $\mathrm{GL}$}
Let $V = V_0 \oplus V_1$ be a super vector space. 
Let $ \mathrm{GL}_V^{sup}$ be the functor which associates to each $R \in 
(\mathsf{Super~Com~Algebras})$ the group of 
$\mathrm{Aut}_R^{sup}(V \otimes R)$ of all super $R$-linear automorphisms on $V \otimes R$.
A \emph{representation} of a super affine group $G$ (or a \emph{super left} $G$-\emph{module} structure) on $V$ 
is a morphism of group-valued
functors $G \to \mathrm{GL}_V^{sup}$. Those representations (or structures) are in a natural one-to-one
correspondence with the super right $\mathcal{O}(G)$-comodule structures on $V$.

\begin{example}
Suppose that $V$ is finite-dimensional, and $m = \dim V_0$, $n = \dim V_1$. Then 
$ \mathrm{GL}_V^{sup}$ is denoted by $\mathrm{GL}(m|n)$. This
is a super algebraic group represented by
$$ \mathcal{O}(\mathrm{GL}(m|n)) = \Bbbk[x_{ij}, y_{k\ell}, \mathrm{det}(X)^{-1},\mathrm{det}(Y)^{-1}] 
\otimes \wedge(p_{i\ell}, q_{kj}). $$
Here $x_{ij}$, $y_{k\ell}$ are even, and $p_{i\ell}$, $q_{kj}$ odd; we suppose that they are entries of the 
matrix 
$$\left( \begin{matrix} X & P \\ Q & Y \end{matrix} \right)=
\left( \begin{matrix} x_{ij} & p_{i\ell} \\ q_{kj} & y_{k\ell} \end{matrix} \right),\ 1\le i, j \le m,\ 1 \le k, \ell \le n.$$
$\wedge(p_{i\ell}, q_{kj})$ denotes
the exterior algebra on the vector space with basis $p_{i\ell}$, $q_{kj}$. 
We choose bases $v_1,\dots, v_m$ of $V_0$ and $v_{m+1},\dots, v_{m+n}$ of $V_1$. Note that every automorphism 
$\sigma \in \mathrm{Aut}_R^{sup}(V \otimes R)$ then arises uniquely from a super algebra map 
$\gamma: \mathcal{O}(\mathrm{GL}(m|n)) \to R$ so that 
\[  \begin{pmatrix} \sigma v_1\ \dots \ \sigma v_{m+n}\end{pmatrix}\otimes 1 =\begin{pmatrix} v_1\ \dots \ v_{m+n}\end{pmatrix} \otimes
 \left( \begin{matrix} \gamma X & \gamma P \\ \gamma Q & \gamma Y \end{matrix} \right). \]
Just as for the ordinary $\mathrm{GL}$, the coalgebra structure maps are given by
$$ \Delta\left( \begin{matrix} X & P \\ Q & Y \end{matrix} \right)=
\left( \begin{matrix} X & P \\ Q & Y \end{matrix} \right)\otimes \left( \begin{matrix} X & P \\ Q & 
Y \end{matrix} \right),\quad
\varepsilon \left( \begin{matrix} X & P \\ Q & Y \end{matrix} \right) = \left( \begin{matrix} I & O \\ O 
& I \end{matrix} \right). $$
From the equation $\begin{pmatrix} X & P \\ Q & Y \end{pmatrix}
\begin{pmatrix} S(X) & S(P) \\ S(Q) & S(Y) \end{pmatrix}=\begin{pmatrix} I & O \\ O & I \end{pmatrix}$
one sees that the antipode $S$ must be given by
\begin{align*}
S(X)~&= (X-PY^{-1}Q)^{-1},  &S(Y)~&= (Y-QX^{-1}P)^{-1},\\
S(P)~&=-X^{-1}PS(Y),   &S(Q)~&=-Y^{-1}QS(X). 
\end{align*}
One sees the algebraic group associated with $\mathrm{GL}(m|n)$ is
$\mathrm{GL}_{V_0}\times \mathrm{GL}_{V_1}$.
\end{example}

Just as for algebraic groups every super algebraic group can be embedded into some $\mathrm{GL}(m|n)$
as its closed super subgroup. 
 
\subsection{Tensor product decomposition theorem}
To state this key result, let $G =\mathrm{SSp}\, A$ be a super affine group. Then we have the associated
affine group $G_{ev}= \mathrm{Sp}\, \overline{A}$, where $\overline{A}=A/(A_1)$. The cotangent super vector
space $T_{\varepsilon}^*(G)$ of $G$ at $1$ is given by $A^+/(A^+)^2$, where 
$A^+ := \mathrm{Ker}(\varepsilon : A \to \Bbbk)$. One sees that the odd component of $T_{\varepsilon}^*(G)$ 
equals
\begin{equation}\label{eq:W}
W^A:= A_1/A_0^+A_1,
\end{equation}
where $A_0^+= A_0 \cap A^+$. We have the tensor product $\overline{A}\otimes \wedge(W^A)$ of
two super Hopf algebras; see Example \ref{ex:exterior}.  
Forgetting some of the structures we regard this as a left $\overline{A}$-comodule
super algebra with counit. Regard $A$ as such an object; the left $A$-comodule $A$ is then
regarded as a left $\overline{A}$-comodule along the quotient map $A \to \overline{A}$.

\begin{theorem}[Tensor product decomposition~\cite{M1}]\label{thm:decomposition}
There is a counit-preserving isomorphism 
$A \overset{\simeq}{\longrightarrow} \overline{A}\otimes \wedge(W^A)$
of left $\overline{A}$-comodule super algebras. 
\end{theorem}

Isomorphisms such as above are not canonical in general. 
The theorem is basic, and is indeed used to prove most of
the results which will be cited in what follows. 
To prove the theorem, \emph{Hopf crossed products}, 
a natural generalization of crossed products of algebras by groups,
are used; see \cite{M1, M3}.  

Decompositions as above might not have been familiar to super geometers, as could be guessed from
the following example. 

\begin{example}
Recall from the previous example
$$ \mathcal{O}(\mathrm{GL}_V^{sup}) = \Bbbk[x_{ij}, y_{k\ell}, \mathrm{det}(X)^{-1},\mathrm{det}(Y)^{-1}] 
\otimes \wedge(p_{i\ell}, q_{kj}). $$
This does NOT give such a decomposition as above. Replacing $p_{i\ell}$, $q_{kj}$ with the entries in 
\[  ( p'_{i\ell} ) := X^{-1}P,\quad ( q'_{kj} ) := Y^{-1}Q,
\footnote{These matrices are obtained by the Hopf-Module Theorem [M.~E.~Sweedler, \emph{Hopf algebras},
Benjamin, 1969, Theorem 4.1.1, Page 84]. In fact, Theorem \ref{thm:decomposition} tells us that
there exists a left $\overline{A}$-comodule algebra section $i : \overline{A} \to A$ of the 
quotient map $q:A \to \overline{A}$, through which $A$ is regarded as a left $\overline{A}$-module and,
moreover, as a \emph{left} $\overline{A}$-Hopf module.
By the Hopf-Module Theorem the left $\overline{A}$-coinvariants in $A$
are given by $\sum_{(a)}i(\overline{S}(q({a}_{(1)}))a_{(2)}$, $a \in A$, where $\overline{S}$ is the 
antipode of $\overline{A}$. For the Example, chosen as $i$ is
the obvious one $X \mapsto X$, $Y \mapsto Y$. It will remain to see $\overline{S}(X)=X^{-1}$, $\overline{S}(Y)=Y^{-1}$.}\]
one has a decomposition as above, 
\[ \mathcal{O}(\mathsf{GL}_V^{sup}) = \Bbbk[x_{ij}, y_{k\ell}, \mathrm{det}(X)^{-1},\mathrm{det}(Y)^{-1}] 
\otimes \wedge(p'_{i\ell}, q'_{kj}).  \]
\end{example}

\subsection{Faithful flatness}
To give an immediate consequence of Theorem \ref{thm:decomposition}, let $f : A \to B$
be a map of super commutative Hopf algebras. We remark that isomorphisms as given by the theorem can be chosen
so as compatible with $f$, so that we have the commutative diagram
\begin{equation*}
\begin{xy}
(60,0)   *++{A}                                  ="lu",
(90,0)  *++{\overline{A} \otimes \wedge(W^A)}    ="ru",
(60,-15) *++{B}                                          ="ld",
(90,-15)*++{\overline{B} \otimes \wedge(W^B).}             ="rd",
{"lu" \SelectTips{cm}{} \ar @{->}^{\simeq\qquad} "ru"},
{"ld" \SelectTips{cm}{} \ar @{->}^{\simeq\qquad} "rd"},
{"lu" \SelectTips{cm}{} \ar @{->}_{f}               "ld"},
{"ru" \SelectTips{cm}{} \ar @{->}^{\overline{f}\otimes\wedge(W^f)}    "rd"},
{"lu" \ar @{}|{\circlearrowleft} "rd"}
\end{xy}
\end{equation*}
Here note that since the constructions of $\overline{A}$, $W^A$ are functorial, we have maps
$\overline{f} : \overline{A}\to \overline{B}$,\ $W^f : W^A \to W^B$. Now, suppose that $f$
is an inclusion $A \subset B$. Then the commutative diagram shows that the last two maps are injections. 
By the classical result cited in Section \ref{subsec:affine_group}, $\overline{B}$ is faithfully flat over
$\overline{A}$. In addition, $\wedge(W^B)$ is free over $\wedge(W^A)$ on both sides. It follows
that $B$ is faithfully flat over $A$ on both sides.

\subsection{Schemes and sheaves in the super situation}\label{subsec:super_scheme}
The definitions of schemes and (dur) sheaves given in the second paragraph of Section \ref{subsec:G_vs_F} 
are directly generalized
to the super situation; see \cite{Z, MZ}.  
The generalized notion of schemes is named \emph{super schemes}. Dur sheaves and sheaves
are generalized by the notions with the same names. 
To generalize the former the exact diagram \eqref{eq:fpqf_covering} should be replaced by the one 
that arises from a morphism $R \to T$ in $(\mathsf{Super~Com~Algebras})$ such that $T$ satisfies those
three equivalent conditions for faithful flatness over $R$ which are given in Lemma \ref{lem:faithfully_flat};
the conditions are now equivalent to the opposite-sided variants due to the super-commutativity assumption. 
Super schemes are necessarily dur sheaves, which in turn are sheaves; see \cite{MZ}.   

Suppose that $G$ is a super affine group. 
Just as in the non-super situation, the faithful flatness result obtained in the last subsection
justifies us to define a \emph{quotient super group} of $G$ to be  
a super affine group which is represented by a super Hopf subalgebra of $\mathcal{O}(G)$. 
We can prove that there is a natural one-to-one correspondence between 
the closed normal super subgroups $N$ of $G$ and the quotient super groups of $G$; 
the quotient super group corresponding to an $N$ is given by the \emph{dur sheafification} $G\tilde{\tilde{/}} N$
of the functor which associates to $R \in (\mathsf{Super~Com~Algebras})$ the quotient group $G(R)/N(R)$.
If $G$ is a super algebraic group, then $\mathcal{O}(G\tilde{\tilde{/}} N) \subset \mathcal{O}(G)$ 
is an fppf covering, so that
$G\tilde{\tilde{/}} N$ coincides with the sheafification $G\tilde{/} N$. 
See \cite{M1, Z}. 

The geometric viewpoint defines a \emph{super space} (over $\Bbbk$) to be a pair $(X, \mathcal{O}_X)$
of a topological space $X$ and a sheaf $\mathcal{O}_X$ of super commutative algebra on $X$ such that
each stalk $\mathcal{O}_{X,x}$ is local, i.e. its even component is local. A \emph{super scheme} is then defined
to be a super space which has an open covering of affine super spaces; see Manin \cite[Chapter~4]{Manin}.  
The Comparison Theorem cited in Section \ref{subsec:G_vs_F} is generalized to the super situation, so that
the two notions of super schemes defined separately from geometric and functorial viewpoints are equivalent;
see \cite[Theorem~5.14]{MZ}. 

\section{The quotient sheaf $G \tilde{/} H$}

\emph{Whether can one generalize Theorem \ref{thm:G/H} to the super situation?} 
This is a question which was posed by Brundan to Zubkov, privately. 
We will answer this in the positive. 

Let $G$ be a super algebraic group, and let $H$ be a closed super subgroup of $G$. 

\begin{theorem}[\cite{MZ}]\label{thm:superG/H}
(1)\ The sheafification $G\tilde{/} H$ of the functor defined on $(\mathsf{Super~Com~Algebras})$
which associates to each $R$ the set $G(R)/H(R)$ of left cosets is a Noetherian super scheme
such that the natural epimorphism $G \to G\tilde{/} H$ is affine and faithfully flat. 

(2)\ The functor $(G\tilde{/} H)_{ev}$ defined on $(\mathsf{Com~Algebras})$ which is obtained from $G\tilde{/} H$
by restriction of the domain is naturally isomorphic to the scheme $G_{ev} \tilde{/} H_{ev}$. 
\end{theorem}

Given $G$, $H$ as above, Brundan \cite{Br} assumes the existence of a Noetherian super scheme $X$ with 
the same properties as those of $G\tilde{/} H$ which are just shown by Part 1 above, 
in order to discuss sheaf cohomologies $H^i(X, \ )$. The assumption was superfluous.

Brundan also asked to Zubkov whether $G\tilde{/} H$ is affine (or representable) whenever the
algebraic group $H_{ev}$ is geometrically reductive. We answer this again in the positive in a
generalized form, as follows; it is known that under the assumption, the scheme $G_{ev} \tilde{/} H_{ev}$ is
affine. 

\begin{prop}[\cite{MZ}]\label{prop:affineG/H}
$G\tilde{/} H$ is affine if and only if $G_{ev} \tilde{/} H_{ev}$ is. 
\end{prop}

Just as in the non-super situation, $G\tilde{/} H$ is affine if and only if $\mathcal{O}(G)$
is an injective cogenerator (or equivalently, faithfully coflat) as a left or right $\mathcal{O}(H)$-comodule; see 
\cite{Z}.  

As was just seen the affinity of $G\tilde{/}H$ is reducible to the same property
of the associated non-super object. 
Three more such properties will be seen 
in Sections \ref{subsec:SC}, \ref{subsec:integral} and \ref{subsec:unipotency}.

\section{Super hyperalgebras and super Lie algebras }

Let us see what roles these super objects play.

\subsection{Super hyperalgebras}
Let $G$ be a super algebraic group. The \emph{super hyperalgebra} $\mathrm{hy}(G)$
of $G$ is defined by the same formula as \eqref{eq:hyperalgebra} when we suppose $A =\mathcal{O}(G)$. 
This is now a super cocommutative
Hopf algebra which is irreducible as a coalgebra. The super vector subspace of $\mathrm{hy}(G)$ consisting
of all primitives
\[ \mathfrak{g} =\{ u \in \mathrm{hy}(G)\mid \Delta(u) = 1 \otimes u + u \otimes 1 \} \]
turns into a finite-dimensional super Lie algebra with respect to 
$[\ , \ ] := \mathrm{product}\circ
(\mathrm{id}_{\mathfrak{g}^{\otimes 2}} - c_{\mathfrak{g}, \mathfrak{g}})$. This 
is called the \emph{super Lie algebra} of $G$, denoted by $\mathrm{Lie}(G)$. We have
\begin{equation}\label{eq:Lie_0}
\mathrm{Lie}(G)_0 = \mathrm{Lie}(G_{ev}).  
\end{equation}
If $\mathrm{char}\, \Bbbk = 0$,
then $\mathrm{hy}(G)$ coincides with the universal envelope $U(\mathrm{Lie}(G))$ of $\mathrm{Lie}(G)$.

The canonical pairing $\mathrm{hy}(G)\times \mathcal{O}(G)
\to \Bbbk$ induces a natural map
\begin{equation}\label{eq:cano}
\mathcal{O}(G) \to \mathrm{hy}(G)^{\circ}
\end{equation}
of super Hopf algebras. Here we remark that given a super Hopf algebra $B$, the dual coalgebra
$B^{\circ}$ defined by \eqref{eq:dual_coalgebra} is naturally a super Hopf algebra. 

Set $W := W^A$ with $A = \mathcal{O}(G)$ (see \eqref{eq:W}), and choose an isomorphism 
$\mathcal{O}(G) \overset{\simeq}{\longrightarrow} \mathcal{O}(G_{ev}) \otimes \wedge(W)$
such as given by Theorem \ref{thm:decomposition}. This induces a unit-preserving isomorphism 
$\mathrm{hy}(G_{ev})\otimes \wedge(W)^* \overset{\simeq}{\longrightarrow}\mathrm{hy}(G)$
of left $\mathrm{hy}(G_{ev})$-module super coalgebras. 
The natural inclusion $\mathrm{hy}(G_{ev}) \subset \mathrm{hy}(G)$ induces a super Hopf
algebra map (indeed, surjection)  $\mathrm{hy}(G)^{\circ} \to \mathrm{hy}(G_{ev})^{\circ}$,
along which $\mathrm{hy}(G)^{\circ}$ will be regarded as a left $\mathrm{hy}(G_{ev})^{\circ}$-comodule.
One can show that the last isomorphism induces a counit-preserving,  
left $\mathrm{hy}(G_{ev})^{\circ}$-comodule super algebra isomorphism 
$\mathrm{hy}(G)^{\circ} \overset{\simeq}{\longrightarrow} \mathrm{hy}(G_{ev})^{\circ}\otimes \wedge(W)$
which fits into the commutative diagram 
\begin{equation*}
\begin{xy}
(0,0)   *++{\mathsf{hy}(G)^{\circ}}                                  ="lu",
(30,0)  *++{\mathsf{hy}(G_{ev})^{\circ}\otimes\wedge(W)}    ="ru",
(0,-15) *++{\mathcal{O}(G)}                                          ="ld",
(30,-15)*++{\mathcal{O}(G_{ev})\otimes\wedge(W).}             ="rd",
{"lu" \SelectTips{cm}{} \ar @{->}^{\simeq\qquad} "ru"},
{"ld" \SelectTips{cm}{} \ar @{->}^{\simeq\qquad} "rd"},
{"ld" \SelectTips{cm}{} \ar @{->}                "lu"},
{"rd" \SelectTips{cm}{} \ar @{->}                "ru"},
{"lu" \ar @{}|{\circlearrowleft} "rd"},
\end{xy}
\end{equation*}
Here the right vertical arrow denotes the natural Hopf algebra map $\mathcal{O}(G_{ev}) \to
\mathrm{hy}(G_{ev})^{\circ}$ tensored with the identity map on $\wedge(W)$. 

A \emph{finite etale group} is an algebraic group such that the corresponding Hopf algebra is finite-dimensional 
separable as an algebra. The super algebraic group $G$ has the largest finite etale quotient group denoted by
$\pi_0(G)$. The following are equivalent:
\begin{itemize}
\item[(i)] $\pi_0(G)$ is trivial;
\item[(ii)] The associated algebraic group $G_{ev}$ is connected;
\item[(iii)] The prime spectrum $\mathrm{Spec}(A_0)$ of the even component
of $A=\mathcal{O}(G)$ is connected;
\item[(iv)] The natural map $\mathcal{O}(G) \to \mathrm{hy}(G)^{\circ}$ given above is injective.
\end{itemize}
If these are satisfied we say that $G$ is \emph{connected}.

Assume that $G$ is connected. Then we may regard $\mathcal{O}(G)\subset \mathrm{hy}(G)^{\circ}$,\ 
$\mathcal{O}(G_{ev})\subset \mathrm{hy}(G_{ev})^{\circ}$ as super Hopf subalgebras,  
via the natural maps. The last commutative
diagram shows the following.

\begin{lemma}[\cite{MS}]\label{lem:characterization}
If $G$ is connected, then $\mathcal{O}(G)$ is characterized in the left $\mathrm{hy}(G_{ev})^{\circ}$-comodule
$\mathrm{hy}(G)^{\circ}$ as the largest $\mathcal{O}(G_{ev})$-subcomodule. 
\end{lemma}

Assume that in addition, $G_{ev}$ is a reductive algebraic group with a split maximal torus $T$. 
Then we have the purely even super Hopf subalgebra $\mathrm{hy}(T)$ of $\mathrm{hy}(G)$. 
The next result follows from Lemma \ref{lem:characterization} combined with the corresponding
result \cite[Part~II, 1.20]{Jantzen} in the non-super situation; see also Section \ref{subsec:dual_coalgebra}. 

\begin{prop}[\cite{MS}]\label{prop:hyG-module}
Given a super vector space $V$, there is a natural one-to-one correspondence between 
\begin{itemize}
\item the super $G$-module structures on $V$, and
\item those locally finite super $\mathrm{hy}(G)$-module structures on $V$ whose 
restricted $\mathrm{hy}(T)$-module structures arise (uniquely) from $T$-module structures.
\end{itemize}
\end{prop}

The last condition on the restricted $\mathrm{hy}(T)$-module structures means that $V$ decomposes
so as $V = \bigoplus_{\lambda \in X(T)}V_{\lambda}$ into weight spaces $V_{\lambda}$, where $X(T)$ denotes
the character group of $T$. 
The result above was previously 
known only for some special super algebraic groups that satisfy the assumption;
see \cite{BrKl, BrKuj, ShuWang}.  

\begin{rem}
The definitions of (super) affine/algebraic groups make sense over any commutative ring. Theory of algebraic
groups over a commutative ring has been established. Indeed, the result \cite[Part~II, 1.20]{Jantzen} cited above
is formulated so as to hold over any integral domain. Accordingly,
Proposition \ref{prop:hyG-module} can be re-formulated in the same situation; see \cite{MS}. 
\end{rem}

\subsection{Harish-Chandra pairs}
Given a super affine group $G = \mathrm{SSp}\, A$, 
the tensor product decomposition theorem tells us that $A$ can recover from $\overline{A}$ and $W^A$
to some extent. One may expect that $A$ or $G$ can recover completely from these two together
with some additional data. This is true if $G$ is a super algebraic group, as will be seen below. 

\begin{definition}[\cite{Kostant, Koszul}]
A \emph{Harish-Chandra pair} 
is a pair $(F, V)$ of an algebraic group $F$ and a
finite-dimensional right $F$-module $V$, given an $F$-module map $[\ , \ ] : V \otimes V \to
\mathrm{Lie}(F)$ such that
\begin{itemize}
\item[(i)] $[u, v]=[v,u]$ for all $u, v \in V$, and
\item[(ii)] $v \triangleleft [v,v]$ for all $v \in V$.
\end{itemize} 
Here, $\mathrm{Lie}(F)$ is regarded as a right $F$-module by the right adjoint action 
(which arises from the conjugation by $F$ on itself), and the $\triangleleft$ in (ii) 
denotes the action
by $\mathrm{Lie}(F)$ on $V$ which arises from the given $F$-module structure on $V$. 
\end{definition}

The Harish-Chandra pairs naturally from a category $(\mathsf{Harish\text{-}Chandra~Pairs})$. 

Let $(\mathsf{Super~Algebraic~Groups})$ denote the category of super algebraic groups, and choose
an object $G = \mathrm{SSp}\, A$ from it. Let $V_G = \mathrm{Lie}(G)_1$ be the odd component of 
$\mathrm{Lie}(G)$; note that $V_G = (W^A)^*$.
The right adjoint action by $G$ on $\mathrm{Lie}(G)$, restricted to $G_{ev}$, stabilizes $V_G$, so that
we have a right $G_{ev}$-module $V_G$. Restrict the bracket on $\mathrm{Lie}(G)$ onto $V_G \otimes V_G$
to obtain 
$[\ , \ ] : V_G \otimes V_G \to \mathrm{Lie}(G_{ev})$. Then the pair 
$(G_{ev}, V_G)$ together with $[\ , \ ]$ just obtained is a Harish-Chandra pair. 

\begin{theorem}\label{thm:HCpair}
$G \mapsto (G_{ev}, V_G)$ gives a category equivalence $$(\mathsf{Super~Algebraic~Groups})\approx 
(\mathsf{Harish\text{-}Chandra~Pairs}).$$
\end{theorem}  

This is a reformulation of \cite[Theorem~29]{M2}, and has the 
same formulation as the corresponding 
result by Koszul \cite{Koszul} (see also \cite[Section~7.4]{Carmeli-Caston-Fioresi}) 
for super Lie groups; the theorem was previously proved by Carmeli and Fioresi \cite{CF} for 
super algebraic groups over an algebraically closed of characteristic zero. 

Let us construct a quasi-inverse of the functor above by using Hopf-algebraic techniques. 
Let $(F, V)$ be a Harish-Chandra pair. Then $V$ is a right $F$-module, whence it is a right Lie module
over the Lie algebra $\mathfrak{g}_0:= \mathrm{Lie}(F)$. There is associated the super Lie 
algebra $\mathfrak{g}_0 \ltimes V$ of semi-direct sum, with even component 
$\mathfrak{g}_0$ and odd component $V$. Note that the bracket on $\mathfrak{g}_0 \ltimes V$ restricted
to $V\otimes V$ is constantly zero. 
Replace this zero map with the
$[\ , \ ]$ associated with the Harish-Chandra pair. Obtained is a new super Lie algebra, say $\mathfrak{g}$.
Let $U_0:= U(\mathfrak{g}_0)$. 
Note that the right $U_0$-module structure on $V$ uniquely gives rise to a right 
$U_0$-module super Hopf algebra structure on the tensor algebra $T(V)$ on $V$, in
which every element in $V$ is supposed to be odd primitive. The associated semi-direct (or smash) product
\[ \mathcal{H}:=U_0 \ltimes T(V) \]
is a super cocommutative Hopf algebra; this is the tensor product $U_0 \otimes T(V)$ as a super coalgebra. 
Note that $U(\mathfrak{g})$ is constructed as the quotient super Hopf algebra of $\mathcal{H}$
divided by the super Hopf ideal
\[ \mathcal{I} = (uv + vu - [u,v]\mid u, v \in V)\]
generated by the indicated even primitives. 

We are going to dualize this last construction. Let $C:=\mathcal{O}(F)$.
Note that the dual space $V^*$ of $V$ is a left $F$-module, or a right $C$-comodule.
The right $C$-comodule structure on $V^*$ uniquely 
gives rise to a right $C$-comodule super Hopf algebra structure
on the graded dual $T_c(V^*):= \bigoplus_{n=0}^{\infty}T^n(V)^*$ of $T(V)$. The associated
smash coproduct $\mathcal{A}:=C \lbtimess T_c(V^*)$ 
is a super commutative Hopf algebra; this is the tensor product $C \otimes T_c(V^*)$ as a super
algebra, and the coproduct on $\mathcal{A}$ is the left $C$-module map, generalizing
\eqref{eq:smash_coproduct}, defined by
\[ \Delta(1 \otimes x) = \sum_{(x)} 1 \otimes \rho(x_{(1)}) \otimes x_{(2)},\quad x \in T_c(V^*), \] 
where $\rho : T_c(V^*)\to T_c(V^*)\otimes C$ denotes the $C$-comodule structure on $T_c(V^*)$. 
This $\mathcal{A}$ is completed to
\[  \widehat{\mathcal{A}} : = (C \lbtimess T_c(V^*))^{\wedge}=\prod_{n=0}^{\infty}C \otimes T^n(V)^* \]
with respect to the linear topology on $\mathcal{A}$ naturally 
given by the $\mathbb{N}$-grading. This $\widehat{\mathcal{A}}$ is a complete topological super commutative 
Hopf algebra, and is a left $C$-comodule super algebra along the projection 
$\widehat{\mathcal{A}} \to C$ onto the zero-th component.
Let $\lambda : \widehat{\mathcal{A}}\to C \otimes \widehat{\mathcal{A}}$ denote the structure map. 
The tensor product of the natural pairings on $C \times U_0$ and on $ T_c(V^*)\times T(V)$ 
gives a pairing $\langle \ , \ \rangle : \widehat{\mathcal{A}} \times \mathcal{H} \to \Bbbk$. 
Being a left $F$-module, $C$ is a left $U_0$-module. Let 
$\mathrm{Hom}_{U_0}(\mathcal{H}, C)$ denote the vector space of all left $U_0$-module
maps $\mathcal{H} \to C$. One sees easily that a linear isomorphism 
\[ \xi : \widehat{\mathcal{A}} \overset{\simeq}{\longrightarrow} \mathrm{Hom}_{U_0}(\mathcal{H}, C)\]
is given by $\xi(a)(x) = \sum_ic_i\langle a_i, x \rangle$, where $a \in \widehat{\mathcal{A}}$,\ 
$x \in \mathcal{H}$ and $\lambda(a) = \sum_ic_i\otimes a_i$.  Transfer the structures on 
$\widehat{\mathcal{A}}$ onto $\mathrm{Hom}_{U_0}(\mathcal{H}, C)$ via $\xi$. One can describe explicitly
the transferred structures, and sees that 
\[ \mathrm{Hom}_{U(\mathfrak{g}_0)}(U(\mathfrak{g}) , C) = \mathrm{Hom}_{U_0}(\mathcal{H}/\mathcal{I}, C)\]
is a discrete super Hopf subalgebra of $\mathrm{Hom}_{U_0}(\mathcal{H}, C)$, and is indeed a super affine
Hopf algebra. The association of the corresponding super algebraic group to $(F, V)$ gives the desired
quasi-inverse. 

\begin{rem}
(1)\ Theorem~\ref{thm:HCpair} shows us a systematic 
method to construct super algebraic groups, and it is applied to
prove Propositions \ref{prop:SC}, \ref{prop:linearly_reductive_positive_char} and Theorem \ref{thm:unipotency}
below. With some modification
the theorem can be re-formulated so as to hold over any commutative ring; see Gavarini \cite[Theorem~4.3.14]{G3}.  

(2)\ Fioresi and Gavarini \cite{FG, G1, G2} constructed \emph{super Chevalley groups} over $\mathbb{Z}$ from 
simple super Lie algebras over $\mathbb{C}$ and their faithful representations; 
they are $\mathbb{Z}$-forms of an important
class of super algebraic groups over $\mathbb{C}$. The re-formulated Theorem~\ref{thm:HCpair} 
gives an alternative, hopefully more conceptual construction of the super Chevalley groups over $\mathbb{Z}$;
see \cite[Sections 9, 10]{MS}. 
\end{rem}

\subsection{Simply-connectedness}\label{subsec:SC}
Given a connected super algebraic group $G$, 
an \emph{etale covering} of $G$ is a pair $(\widetilde{G}, \eta)$ of a connected super algebraic group $\widetilde{G}$ and an
epimorphism $\eta : \widetilde{G} \to G$ of super algebraic groups such that the kernel $\mathrm{Ker}\, \eta$ is finite etale.
We say that $G$ is \emph{simply connected} if it has no non-trivial etale covering.

\begin{prop}[\cite{M2}]\label{prop:SC} 
A connected super algebraic group $G$ is simply connected if and only if $G_{ev}$ is simply connected.
\end{prop}

Suppose that (a)~$\Bbbk$ is an algebraically closed field of characteristic zero, or
(b)~$\Bbbk$ is a perfect field of characteristic~$>2$. It follows from Proposition \ref{prop:SC} 
combined results by Hochschild \cite{H} (in case (a)) and by Takeuchi \cite{T2, T3} (in case (b)) that for a simply
connected super algebraic group $G$, $\mathcal{O}(G)$ is described in terms of 
$\mathfrak{g}:=\mathrm{Lie}(G)$ (in case (a)) or of $\mathrm{hy}(G)$ (in case (b)).
The result is simpler in case (b), and is then that $\mathcal{O}(G) = \mathrm{hy}(G)^{\circ}$. 
In case (a),   
$\mathcal{O}(G)$ is the super Hopf subalgebra of $U(\mathfrak{g})^{\circ}$
consisting of those elements each of which annihilates an ideal of $U(\mathfrak{g})$
generated by some power $\mathrm{Rad}(\mathfrak{g}_0)^n$, $0< n \in \mathbb{Z}$, of  the radical
$\mathrm{Rad}(\mathfrak{g}_0)$ of $\mathfrak{g}_0$. Therefore, if $G_{ev}$ is supposed to be
semisimple in addition, then we have $\mathcal{O}(G) = U(\mathfrak{g})^{\circ}$. 

\section{Linear reductivity and unipotency}

We will characterize super affine or algebraic groups with representation theoretic properties
such as above. 

\subsection{Linear reductivity in characteristic zero}\label{sec:linearly_reductive} 
A super affine group $G$ is said to be \emph{linearly reductive} if every 
super $G$-module is semisimple. 
Every quotient super group of a linearly reductive
affine super group is linearly reductive. 

Suppose that $\mathrm{char}\, \Bbbk = 0$. In the non-super situation 
the Chevalley Decomposition Theorem states
that any affine group $G$ is a semi-direct product $G_r \ltimes G_u$ of the unipotent
radical $G_u$ by a linearly reductive affine group $G_r$. Therefore, if $G_u$ is trivial 
(especially, if $G$ is a reductive algebraic group), then $G$ is linearly reductive. Contrary to
this, linearly reductive super algebraic groups are rather restricted, as will be seen below.

Note that every super affine group $G = \mathrm{SSp}\, A$ has the largest purely
even quotient super group $G_{qev}$; one sees that $\mathcal{O}(G_{qev})$ is the pull-back
$\Delta^{-1}(A_0\otimes A_0)$ of $A_0\otimes A_0$ along the coproduct on $A$. If $G$ is a
super algebraic group, the largest finite etale quotient $G \to \pi_0(G)$ factors through $G_{qev}$.

We say that a super algebraic group $G$ is \emph{tight} if the adjoint action by $G_{ev}$ on the
odd component of $\mathrm{Lie}(G)$ is faithful; the condition is equivalent to that $\mathcal{O}(G)$
is the smallest super Hopf subalgebra of $\mathcal{O}(G)$ that includes the $G_{ev}$-invariants 
$\mathcal{O}(G)^{G_{ev}}$ in $\mathcal{O}(G)$. One sees that every super algebraic 
group $G$ has the largest tight quotient super group, which we denote by $G_{ti}$.

The following is a reformulation of Weissauer's Theorem.

\begin{theorem}[Weissauer~\cite{W}]\label{thm:Weissauer}
Assume that $\Bbbk$ is an algebraically closed field of characteristic zero.
\begin{itemize}
\item[(1)] Those linearly reductive super algebraic group which are tight and connected
are exhausted by finite products
\[ \prod_{r}Spo(1, 2r)^{n_r}, \quad n_r \ge 0 \]
of the ortho-symplectic super algebraic groups $Spo(1, 2r)$, $r>0$. 
\item[(2)] Those linearly reductive super algebraic group which are tight
are exhausted by semi-direct products
\[ \Gamma \ltimes \prod_{r}Spo(1, 2r)^{n_r}, \]
where $\Gamma$ is a finite group, and  
acts on each product $Spo(1, 2r)^{n_r}$ via permutations of components so that
the resulting group map $\Gamma \to \prod_r\mathfrak{S}_{n_r}$ is injective. 
\item[(3)] Every linearly reductive algebraic super group $G$ is naturally isomorphic to the fiber product
\[ G_{qev} \times_{\pi_0(G_{ti})} G_{ti}, \]
where $G_{ti}$ is such as given in (2) with $\Gamma = \pi_0(G_{ti})$. 
\end{itemize}
\end{theorem}

Keep $\Bbbk$ as assumed by the theorem. 
$Spo(1, 2r)$ is the super algebraic group which corresponds to the Harish-Chandra
pair $(Sp_{2r}, V)$ defined by the following. 
\begin{itemize}
\item
$Sp_{2r}$ is the symplectic group\footnote{To be more precise, $Sp_{2r}$ is meant to be the algebraic group which 
arises from the linear algebraic group described here; see the second paragraph of 
Section \ref{subsec:algebraic_group}.} of degree $2r$, which thus consists of the matrices
$g\in \mathrm{GL}_{2r}(\Bbbk)$ such that $g\, J \, {}^tg = J$, where $J$ is a fixed anti-symmetric matrix
in $\mathrm{GL}_{2r}(\Bbbk)$; 
the Lie algebra $sp_{2r}=\mathrm{Lie}(Sp_{2r})$ of $Sp_{2r}$ consists of the $2r\times 2r$ matrices $X$ such that
$XJ$ is symmetric. 
\item
$V$ is the vector space $k^{2r}$ of row vectors with $2r$ entries, and is
regarded as a right $Sp_{2r}$-module by the matrix multiplication. 
\item
The structure $[\ , \ ] : V \otimes V \to
sp_{2r}$ is defined by 
\[ [v,w] = \frac{1}{2}\, J\big( {}^tv\, w + {}^tw\, v \big) , \quad v, w \in V. \]
\end{itemize}
The super algebraic group $Spo(1, 2r)$ is simple (i.e. does not contain any non-trivial closed normal
super subgroup) and simply connected. 

\subsection{Linear reductivity in positive characteristic} 
In positive characteristic the situation is more restrictive, as is seen from the following.

\begin{prop}[\cite{M2}]\label{prop:linearly_reductive_positive_char}
Assume that
$\mathrm{char}\, \Bbbk >2$. 
Then a linearly reductive super affine group $G$ is necessarily purely even, i.e. is
an ordinary affine group. Hence by Nagata's Theorem, $\mathcal{O}(G)\otimes \overline{\Bbbk}$
is a group algebra provided $G$ is algebraic and connected. 
\end{prop}

To discuss here Frobenius morphisms, let $G = \mathrm{SSp}\, A$ be a super affine group over a field 
$\Bbbk$ of characteristic $p > 2$. Regard the super Hopf algebra 
\[ A^{(p)} : = A \otimes \Bbbk^{1/p} \]
over $\Bbbk^{1/p}$ as a super Hopf algebra over $\Bbbk$ via $c \mapsto 1 \otimes c^{1/p}$, $\Bbbk \to 
A\otimes \Bbbk^{1/p} $, 
and let $G^{(p)} = \mathrm{SSp}\, A^{(p)}$ denote the corresponding super affine group over $\Bbbk$. 
One sees that
\[ F_A : A^{(p)} \to A, \quad F_A(a \otimes c)= a^pc^p \]
is a super Hopf algebra map. Let $F_G : G \to G^{(p)}$ denote the corresponding morphism of super
affine groups. This is called the \emph{Frobenius morphism} for $G$, and the kernel $G_p:=\mathrm{Ker}\, F_G$
is called the \emph{Frobenius kernel} of $G$. Recall $W^A$ from \eqref{eq:W}.
Since the image of $F_{\wedge(W^A)}$ equals $\Bbbk$, we
see from Theorem \ref{thm:decomposition} that the image of $F_A$ is a purely even super Hopf subalgebra
of $A$, and it coincides with the image of $F_{\overline{A}}$ with $\overline{A}=\mathcal{O}(G_{ev})$.  
Therefore, $(G_p)_{ev}$ equals the Frobenius kernel $\mathrm{Ker}\, F_{G_{ev}}$ of $G_{ev}$, and 
the $W^B$ of $B:=\mathcal{O}(G_p)$ equals $W^A$. 
It follows that 
if $G$ is a super algebraic group, then $G_p$ is infinitesimal in the sense that $\mathcal{O}(G_p)$ is
finite-dimensional, and is local, i.e. the augmentation ideal $\mathcal{O}(G_p)^+$ is nilpotent. 

In virtue of the situation above one sometimes finds it easier to prove results in positive characteristic than in 
characteristic zero. 

\subsection{Integrals}\label{subsec:integral}
Let $G$ be a super affine group. A \emph{left integral} for $G$ is a (not necessarily super) left 
$\mathcal{O}(G)$-comodule map $\int : \mathcal{O}(G)\to \Bbbk$. Such an $\int$ is necessarily
homogeneous, i.e. $\int \mathcal{O}(G)_0 = 0$ or $\int \mathcal{O}(G)_1 = 0$. Moreover, 
the left integrals for $G$ form a vector space, say $\mathcal{I}_{\ell}(G)$, of dimension $\le 1$. The vector
space $\mathcal{I}_r(G)$ of the \emph{right integrals} for $G$, which are defined in the obvious manner, 
is isomorphic to $\mathcal{I}_{\ell}(G)$
via an isomorphism $\mathcal{I}_{\ell}(G) \overset{\simeq}{\longrightarrow} \mathcal{I}_{r}(G)$ given by 
$\int\mapsto \int\circ S$. See Scheunert and Zhang \cite{SZ}. 

We will say that $G$ \emph{has an integral} if $\mathcal{I}_{\ell}(G)\ne 0$ or equivalently, if $\mathcal{I}_r(G)\ne 0$. 
It is known that $G$ has an integral if $\dim \mathcal{O}(G)< \infty$. 

The following two lemmas follow easily from the corresponding results on ordinary Hopf algebras,  
by using the bosonization; see Section \ref{subsec:bosonization}. 

\begin{lemma}
The following are equivalent:
\begin{itemize}
\item[(1)] $G$ has an integral;
\item[(2)] $\mathcal{O}(G)\oplus \mathcal{O}(G)[1]$ 
is a generator of the category of super $G$-modules;
\item[(3)] Any injective super $G$-modules is necessarily projective;
\item[(4)] The injective hull of every finite-dimensional super $G$-module is
finite-dimensional. 
\end{itemize}
\end{lemma}

\begin{lemma}
A super affine group $G$ is linearly reductive if and only if there exists a one-sided
(necessarily, two-sided) integral $\int : \mathcal{O}(G) \to \Bbbk$ such that $\int 1 \ne 0$.
\end{lemma}

The last equivalent conditions are rarely satisfied
unless $G$ is an ordinary affine group, as is seen from the results of the last two subsections.

The following is an unpublished result by Taiki~Shibata and myself. 
\begin{prop}
Suppose that $G$ is a super algebraic group. Then the following are equivalent:
\begin{itemize}
\item[(a)] $G$ has an integral;
\item[(b)] The associated algebraic group $G_{ev}$ has an integral.
\end{itemize}
\end{prop}

Sullivan \cite{Sullivan2} (see also \cite{Sullivan1})
tells us that if $\mathrm{char}\, \Bbbk = 0$, then Condition (b) is equivalent to 
that $G_{ev}$ is linearly reductive, and that if $\mathrm{char}\, \Bbbk >2$, then the condition
is equivalent to that the reduced algebraic group $F_{red}$ associated with $F$ is a torus, where $F$ denotes
the identity connected component $(G_{ev})^0_{\overline{\Bbbk}}$
of the base extension $(G_{ev})_{\overline{\Bbbk}}$ of $G_{ev}$ to the algebraic closure $\overline{\Bbbk}$ of $\Bbbk$. 

If $\mathrm{char}\, \Bbbk = 0$, we have thus many examples of super algebraic groups which 
are not linearly reductive, but have integrals. 

\subsection{Unipotency}\label{subsec:unipotency}
A property very opposite to linear reductivity is unipotency.
A super affine group $G$ is said to be \emph{unipotent} if simple super $G$-modules 
are exhausted by the trivial super $G$-module $\Bbbk$ which may be purely even or odd,
or in other words, if $\mathcal{O}(G)$ is irreducible as a coalgebra. 
 
The following result is due to Alexandr~N.~Zubkov, and is contained in \cite{M2}, given a simple proof.

\begin{theorem}[Zubkov]\label{thm:unipotency}
A super affine group $G$ is unipotent if and only if the associated affine group $G_{ev}$ 
is unipotent. 
\end{theorem} 

\section*{Acknowledgments}
The work was supported by
Grant-in-Aid for Scientific Research (C) 23540039, Japan Society of the Promotion of Science. 
I thank the organizers of the symposium for giving me an opportunity to present the results.

\end{document}